\documentstyle[12pt]{article}
\newcommand\qed{{\unskip\nobreak\hfil\penalty50\hskip2em\vadjust{}
    \nobreak\hfil$\Box$\parfillskip=0pt\finalhyphendemerits=0\par}}

\renewcommand{\thesection}{\arabic{section}}
\renewcommand{\theequation}{\thesection.\arabic{equation}}

\newtheorem{theorem}{Theorem}[section]

\newtheorem{proposition}{Proposition}[section]
\newtheorem{definition}{Definition}[section]
\newtheorem{remark}{Remark}[section]

\newcommand{\x}{\times}

\renewcommand{\a}{\alpha}
\renewcommand{\b}{\beta}
\renewcommand{\d}{\delta}
\newcommand{\D}{\Delta}
\newcommand{\e}{\varepsilon}
\newcommand{\g}{\gamma}

\renewcommand{\l}{\lambda}

\newcommand{\var}{\varphi}
\newcommand{\s}{\sigma}

\newcommand{\th}{\theta}

\newcommand{\nb}{\mbox{{\bf N}}{}}
\newcommand{\rb}{\mbox{{\bf R}}{}}
\newcommand{\cb}{\mbox{{\bf C}}{}}

\newcommand{\sd}{\rhd\mbox{\hspace{-2ex}}<}

\title{Quasitriangularity and Enveloping Algebras \\
for Inhomogeneous Quantum Groups}
\author{P. Podle\'s\thanks{On leave
{}from Department of Mathematical Methods in Physics, Faculty of Physics,
University of Warsaw, Ho\.{z}a 74, 00-682 Warszawa,
Poland}\ \thanks{This research was supported in
part by NSF grant DMS-9508597 and in part by Polish KBN grant No. 2 P301 020
07} \\
Department of Mathematics \\
University of California \\
Berkeley, CA\ \ 94720, USA}
\date{}

\begin{document}

\maketitle

\begin{abstract}
Coquasitriangular universal ${\cal R}$ matrices on quantum Lorentz and quantum
Poincar\'e groups are classified.  The results extend (under certain
assumptions) to inhomogeneous quantum groups of \cite{INH}.  Enveloping
algebras on those objects are described.
\end{abstract}

\setcounter{section}{-1}
\section{Introduction}
\label{sec0}

Possible $R$-matrices for the fundamental representations of inhomogeneous
quantum groups were found in Proposition 3.14 of
\cite{INH}.  In the present paper we describe
universal ${\cal R}$ matrices for those objects (under certain assumptions
which
are fulfilled in the case of quantum Poincar\'e groups \cite{POI}).  Our study
will be useful in developing a simple physical model \cite{DM} of free
particles on a quantum Minkowski space \cite{POI}.

In Section~\ref{sec1} we show how to construct co(quasi)triangular
(${}^*$-)bialgebras and Hopf\linebreak
 (${}^*$-)algebras whose relations are given by
means of
intertwiners:  we simplify and extend results of \cite{CQG}.  Next, in
Section~\ref{sec2} we classify co(quasi)triangular (${}^*$-)structures
${\cal R}$
on quantum Lorentz groups \cite{WZ}.  Using the results of Sections~\ref{sec1}
and \ref{sec2}, in Section~\ref{sec3} we classify such structures on quantum
Poincar\'e groups \cite{POI} and also, more generally, on inhomogeneous quantum
groups \cite{INH} (certain natural assumptions regarding
mainly restriction of those
structures to the homogeneous quantum group are made).  In Section~\ref{sec4}
we use ${\cal R}$ to define enveloping algebras for inhomogeneous quantum
groups.

We sum over repeated indices (Einstein's convention).  We work over the field
${\cb}$.  We write $a \sim b$ if $a = \l b$ for $\l \in {\cb
}{\backslash}\{0\}$.  If $V,W$ are vector spaces then $\tau: V \otimes W \to W
\otimes V$ is given by $\tau(x\otimes y) = y \otimes x$, $x \in V$, $y \in W$.
If ${\cal A}$ is an algebra, $v \in M_N({\cal A})$, $w \in M_K({\cal A})$, then
the tensor product $v \otimes w \in M_{NK}({\cal A})$ is defined by $(v \otimes
w)_{ij,kl} = v_{ik}w_{jl}$, $i,k = 1,\dots,N$, $j,l = 1,\dots,K$.  We set $\dim
v = N$.  If ${\cal A}$ is a ${}^*$-algebra then the conjugate of $v$
is defined as
${\bar v} \in M_N({\cal A})$ where ${\bar v}_{ij} = v_{ij}{}^*$, $i,j =
1,\dots,N$.

Throughout the paper quantum groups $H$ are abstract objects described by the
corresponding (${}^*$-)bialgebras $\mbox{Poly}(H) = ({\cal A},\D)$.  We
denote by
$\D,\e,S$ the comultiplication, counit and
(if exists) coinverse of $\mbox{Poly}(H)$.  We
say that $v$ is a representation of $H$ (i.e. $v \in \mbox{ Rep } H$) if $v \in
M_N({\cal A})$, $N \in {\nb}$, and $\D v_{ij} = v_{ik} \otimes v_{kj}$,
$\e(v_{ij}) = \d_{ij}$, $i,j = 1,\dots,N$.  The conjugate of a representation
and tensor product of representations are also representations.  Matrix
elements of representations of
$H$ span
${\cal A}$ as a linear space.  The set of nonequivalent irreducible
representations of $H$ is denoted by $\mbox{Irr } H$.  If $v,w \in \mbox{ Rep }
H$, then we say that $A \in M_{\dim w \x \dim v}({\cb})$ intertwines $v$ with
$w$ (i.e. $A \in \mbox{ Mor}(v,w)$) if $Av = wA$.  The dual vector space
${\cal A}'$ is an algebra w.r.t. the convolution defined by $\rho * \rho' =
(\rho
\otimes \rho')\D$.  It has a unit $I_{{\cal A}'} = \e$. For $\rho\in
{\cal A}'$,
$a\in {\cal A}$, we set $\rho * a=( id \otimes\rho)\D a$,
$a * \rho=(\rho\otimes id )\D a$.

\section{Coquasitriangular bialgebras}
\label{sec1}

In this section we discuss bialgebras $\mbox{Poly}(G) = ({\cal B},\D)$
defined by
several fundamental representations of $G$ and intertwiners among them.  We
provide necessary and sufficient conditions for the
existence of coquasitriangular
structure ${\cal R}$ for $G$.  Hopf and $*$ structures are also investigated.
The results generalize
results of Theorem~$1.4$ of \cite{CQG} (see also Theorem~$1.1$ of
\cite{WZ}).

\begin{proposition}
\label{prop1.1}
Let ${\cal B}$ be an algebra generated by $w^{\a}_{mn}$, $m,n =
1,\dots,d_{\a}$, $\a \in {\cal J}$ and relations ($0 \in {\cal J}$, $d_0=1$)
\setcounter{equation}{0}
\renewcommand{\theequation}{\thesection.\arabic{equation}}
\begin{equation}
\label{eq1.1}
w^{\a} \otimes w^0 = w^0 \otimes w^{\a} = w^{\a},\ \a \in {\cal J},
\end{equation}
\begin{equation}
\label{eq1.2}
W_m(w^{\a_{m1}} \otimes w^{\a_{m2}} \otimes \dots w^{\a_{ms_m}}) =
(w^{\b_{m1}} \otimes \dots \otimes w^{\b_{mt_m}})W_m,\ m \in K.
\end{equation}
Then $w^0 = (I_B)$ and there exists a unique $\D$ such that $\mbox{\em
Poly}(G) =
({\cal B},\D)$ is a bialgebra and $w^{\a}$, $\a \in {\cal J}$, representations
of $G$.
\end{proposition}

\medskip
\noindent
{\bf Proof} is the same as in Theorem~$1.4$ of \cite{CQG}. \qed

\medskip
Let us recall

\begin{definition}[cf. \cite{Dr}, \cite{Ha}, \cite{LT}]
\label{def1.2}
{\em We say that $({\cal B},\D,{\cal R})$ is a coquasitriangular (CQT)
bialgebra
if
$({\cal B},\D)$ is a bialgebra and ${\cal R} \in ({\cal B} \otimes {\cal B})'$
satisfies}
\begin{equation}
\label{eq1.3}
{\cal R}(I \otimes x) = {\cal R}(x \otimes I) = \e(x),
\end{equation}
\begin{equation}
\label{eq1.4}
{\cal R}(xy \otimes z) = {\cal R}(x \otimes z^{(1)}) {\cal R}(y \otimes
z^{(2)}),
\end{equation}
\begin{equation}
\label{eq1.5}
{\cal R}(x \otimes yz) = {\cal R}(x^{(1)} \otimes z) {\cal R}(x^{(2)} \otimes
y),
\end{equation}
\begin{equation}
\label{eq1.6}
y^{(1)}x^{(1)}{\cal R}(x^{(2)} \otimes y^{(2)}) = {\cal R}(x^{(1)} \otimes
y^{(1)})x^{(2)}y^{(2)}
\end{equation}
where we have used the Sweedler's notation
$\Delta(x)=x^{(1)}_i\otimes x^{(2)}_i
\equiv x^{(1)}\otimes x^{(2)}$.
\end{definition}

\begin{remark}[cf. \cite{Lyu}, \cite{LT}, Proposition~1.3 of \cite{CQG}]
\label{rem1.3}
{\em Let $\mbox{Poly}(G) = ({\cal B},\D)$ be a bialgebra and ${\cal R} \in
({\cal B}
\otimes {\cal B})'$.  For each $v,w \in \mbox{ Rep } G$ one can
define
$R^{vw} \in
\mbox{ Lin}({\cb}^{\dim v} \otimes {\cb}^{\dim w},{\cb}^{\dim w}
\otimes {\cb}^{\dim v})$ by}
\begin{equation}
\label{eq1.6'}
(R^{vw})_{ij,kl} = {\cal R}(v_{jk} \otimes w_{il}),\ j,k = 1,\dots,\dim v,\ i,l
= 1,\dots,\dim w.
\end{equation}
Then
\begin{equation}
\label{eq1.6a}
(1\!\!1 \otimes S)R^{v_1w} = R^{v_2w}(S \otimes 1\!\!1)
\mbox{ {\em if $S \in$ Mor$(v_1,v_2)$,}}
\end{equation}
\begin{equation}
\label{eq1.6b}
(S \otimes 1\!\!1)R^{vw_1} = R^{vw_2}(1\!\!1 \otimes S)
\mbox{ {\em if $S \in$ Mor$(w_1,w_2)$,}}
\end{equation}
{\em $v,w,v_1,w_1 \in \mbox{ Rep } G$.  Suppose that $L \subset \mbox{ Rep } G$
is  such that
the matrix elements of representations $v \in L$ linearly span $\cal
B$ (e.g. $L = \mbox{ Rep } G$).  Consider the conditions}
\begin{equation}
\label{eq1.7}
R^{0v} = R^{v0} = 1\!\!1,\ v \in L,
\end{equation}
\begin{equation}
\label{eq1.8}
R^{v_1 \otimes v_2,w} = (R^{v_1w} \otimes 1\!\!1)(1\!\!1 \otimes R^{v_2w}),\
v_1,v_2,w \in L,
\end{equation}
\begin{equation}
\label{eq1.9}
R^{v,w_1\otimes w_2} = (1\!\!1 \otimes R^{vw_2})(R^{vw_1} \otimes 1\!\!1),\
v,w_1,w_2 \in L,
\end{equation}
\begin{equation}
\label{eq1.10}
R^{vw} \in \mbox{ {\rm Mor}}(v \otimes w,w \otimes v),\ v,w \in L.
\end{equation}
{\em Then (\ref{eq1.7}) $\Leftrightarrow$ (\ref{eq1.3}), (\ref{eq1.8})
$\Leftrightarrow$ (\ref{eq1.4}), (\ref{eq1.9}) $\Leftrightarrow$ (\ref{eq1.5}),
(\ref{eq1.10}) $\Leftrightarrow$ (\ref{eq1.6})}.
\end{remark}

\begin{theorem}
\label{th1.4}
Let $\mbox{\em Poly}(G) = ({\cal B},\D)$ be a bialgebra defined in
Proposition~\ref{prop1.1} and $R^{\a\b} \in \mbox{ Lin}({\cb}^{d_{\a}}
\otimes {\cb
}^{d_{\b}},{\cb}^{d_{\b}} \otimes {\cb}^{d_{\a}})$, $\a,\b\in {\cal J}$.

The following are equivalent

1)  there exists ${\cal R} \in ({\cal B} \otimes {\cal B})'$ such that
$({\cal B},\D,{\cal R})$ is a CQT bialgebra and $R^{\a\b} = R^{w^{\a}w^{\b}}$.

2)
\begin{equation}
\label{eq1.11}
R^{0\a} = R^{\a 0} = 1\!\!1,\ \a \in {\cal J},
\end{equation}
\begin{equation}
\label{eq1.12}
(1\!\!1 \otimes W_m)R^{\a_{m1}\cdot \dots \cdot \a_{ms_m},\g} =
R^{\b_{m1}\cdot \dots \cdot \b_{mt_m},\g}(W_m \otimes 1\!\!1),\ m \in K,\ \g
\in {\cal J}\setminus\{0\},
\end{equation}
\begin{equation}
\label{eq1.13}
R^{\g,\b_{m1}\cdot \dots \cdot \b_{mt_m}}(1\!\!1 \otimes W_m) = (W_m \otimes
1\!\!1)R^{\g,\a_{m1}\cdot \dots \cdot \a_{ms_m}},\ m \in K,\ \g \in {\cal J}
\setminus\{0\},
\end{equation}
\begin{equation}
\label{eq1.14}
R^{\a\b} \in \mbox{ Mor}(w^{\a} \otimes w^{\b},\ w^{\b} \otimes w^{\a}),\ \a,\b
\in {\cal J}\setminus\{0\},
\end{equation}
where $R^{\d_1\cdot \dots \cdot \d_{k+1},\g} = (R^{\d_1 \cdot \dots \cdot
\d_k,\g} \otimes 1\!\!1)(1\!\!1 \otimes R^{\d_{k+1}\g})$, $R^{\g,\d_1 \cdot
\dots \cdot \d_{k+1}} = (1\!\!1 \otimes R^{\g\d_{k+1}})(R^{\g,\d_1\cdot \dots
\cdot \d_k} \otimes 1\!\!1)$, $\g,\d_1,\dots,\d_{k+1} \in {\cal J}$, $k =
1,2,\dots$.

Moreover, such ${\cal R}$ is unique.
\end{theorem}

\medskip
{\bf Remark.} In special cases related statements can be found in \cite{Ma},
\cite{LT}, \cite{Ha} and Theorem 1.4 of \cite{CQG}. Cf.  \cite{Lyu}.

\medskip
\noindent
{\bf Proof}.  Assume condition 1).  According to Remark~\ref{rem1.3},
(\ref{eq1.7})--(\ref{eq1.10}) follow.  Thus we get (\ref{eq1.11}),
(\ref{eq1.14}).  Using (\ref{eq1.6a})--(\ref{eq1.6b}) for $S = W_m$ (see
(\ref{eq1.2})) and (\ref{eq1.8})--(\ref{eq1.9}), one obtains
(\ref{eq1.12})--(\ref{eq1.13}) and the condition 2) is proved.

Moreover, $R^{w^{\a}w^{\b}} = R^{\a\b}$ uniquely determine $R^{vw}$, where $v,w
\in L_0$,
\[
L_0 = \{\mbox{tensor products of some number of copies of representations
$w^{\a}$, $\a \in {\cal J}$}\}.
\]
Using (\ref{eq1.6'}) and the fact that the
matrix elements of representations from
$L_0$ linearly span $\cal B$,
the uniqueness of ${\cal R}$ follows.  It remains to
prove

2) $\Rightarrow$ 1):  Assume condition 2).
Using (\ref{eq1.11}), we can replace ${\cal J}\setminus\{0\}$ by ${\cal J}$
in (\ref{eq1.12})-(\ref{eq1.14}). We define the
homomorphisms
${\cal R}^{\b}: {\cal B} \to M_{d_{\b}}({\cb})$, $\b \in {\cal J}$, by
\begin{equation}
\label{eq1.15}
[{\cal
R}^{\b}(w^{\a}_{ij})]_{kl} = R^{\a\b}_{ki,jl}, \quad \a\in{\cal J}
\end{equation}
(later on we will have ${\cal R}_{kl}^{\b} = {\cal R}(\cdot \otimes
w_{kl}^{\b})$).  They preserve the relations (\ref{eq1.1})--(\ref{eq1.2}) due
to (\ref{eq1.11})--(\ref{eq1.12}).  Setting $\a = 0$ in (\ref{eq1.15}),
we show that ${\cal
R}^{\b}$ are unital.  Setting $\b = 0$ in (\ref{eq1.15}),
one gets ${\cal R}_{11}^0 =
\e$ (it is true on
the generators $w_{ij}^{\a}$).  Hence ${\cal R}_{11}^0 * {\cal
R}_{kl}^{\b} = {\cal R}_{kl}^{\b} * {\cal R}_{11}^0 = {\cal R}_{kl}^{\b}$.  Let
\[
\begin{array}{rll}
X &= &\{x \in {\cal B}: W_{b_1\dots b_{t_m},a_1\dots a_{s_m}}({\cal
R}_{a_{s_m}c_{s_m}}^{\a_{ms_m}} * \dots * {\cal R}_{a_1c_1}^{\a_{m1}})(x) \\
&= &({\cal R}_{b_{t_m}d_{t_m}}^{\b_{mt_m}} * \dots * {\cal
R}_{b_1d_1}^{\b_{1t_1}})(x)W_{d_1\dots d_{t_m},c_1 \dots c_{s_m}},\
m \in
K\}.
\end{array}
\]
It is straightforward
to show that $X$ is an algebra.  According to (\ref{eq1.13}),
$u_{kl}^{\g} \in X$ $(k,l = 1,\dots,d_{\g},\ \g \in {\cal J})$.  Thus $X =
{\cal B}$.  Hence there exists a linear antihomomorphism $\th: {\cal
B} \to
{\cal B}'$ such that $\th(w^{\b}_{kl}) = {\cal R}_{kl}^{\b}$, $k,l =
1,\dots,d_{\b}$,
$\b \in {\cal J}$ ($\th$ preserves (\ref{eq1.1})--(\ref{eq1.2})).  Setting
${\cal R}(x
\otimes y) = [\th(y)](x)$, $x,y \in {\cal B}$, we obtain
an ${\cal R} \in ({\cal
B}
\otimes {\cal B})'$ which satisfies (\ref{eq1.5}), (\ref{eq1.3}).
Moreover, (\ref{eq1.15}) yields $R^{\a\b} = R^{w^{\a}w^{\b}}$.
Let $Y = \{z
\in {\cal B}:
\forall x,y \in {\cal B}$ ${\cal R}(xy \otimes z) = {\cal R}(x \otimes
z^{(1)}){\cal R}(y \otimes z^{(2)})\}$.  Then $Y$ is an algebra (we use
(\ref{eq1.5})) and
$w^{\b}_{kl} \in Y$.  Hence, $Y = {\cal B}$ and (\ref{eq1.4}) follows.  Thus
(Remark~\ref{rem1.3}) we get (\ref{eq1.7})--(\ref{eq1.9}) for $L = L_0$.  That
and (\ref{eq1.10}) for $v,w \in \{w^{\a}: \a \in {\cal J}\}$ (see
(\ref{eq1.14})) give (\ref{eq1.10}) for $L = L_0$, hence (Remark~\ref{rem1.3})
(\ref{eq1.6})
and the condition 1) is satisfied. \qed

\begin{remark}
\label{rem1.4'}
{\em The
unital antihomomorphism $\th: {\cal B} \to {\cal B}'$ introduced in the above
proof exists for each CQT bialgebra (cf. \cite{Lyu})
and satisfies $(\th \otimes \th)\D  = \D'\th$
where $\D': {\cal B}' \to ({\cal B} \otimes {\cal B})'$ is defined by
$\D'(\var) = \var \circ m$, $\var \in {\cal B}'$, cf.
Appendix of \cite{CQG}}.
\end{remark}

Let us now pass to the Hopf algebra structure.

\begin{proposition}[cf. \cite{TK}, Proof of Theorem $1.4.1$ of \cite{CQG}]
\label{prop1.5}
Let $\mbox{\em Poly}(G) = ({\cal B},\D)$
be a bialgebra and $w,w',w''$ be representations of
$G$.  Suppose there exist $E \in \mbox{ {\em Mor}}(I,w \otimes w')$, $E' \in
\mbox{ {\em Mor}}(w'' \otimes w,I)$ such that $E$ is left nondegenerate, i.e.
$E_{i-} = (E_{ij})_{j=1}^{\dim w'}$, $i = 1,\dots,\dim w$, are linearly
independent
and  $E'$ is right nondegenerate, i.e. $E'_{-k} = (E'_{mk})_{m=1}^{\dim
w''}$, $k = 1,\dots,\dim w$, are linearly independent.  Then $w^{-1}$ exists.
\end{proposition}

\begin{proposition}[cf. \cite{TK}, Proof of Theorem $1.4.1$ of \cite{CQG}]
\label{prop1.6}
Let $({\cal B},\D)$ be the bialgebra defined in Proposition~\ref{prop1.1}.
Suppose
$(u^{\a})^{-1}$ exist, $\a \in {\cal J}$.  Then $({\cal B},\D)$ is a Hopf
algebra.
\end{proposition}

A CQT bialgebra $({\cal B},\D,{\cal R})$ such that $(B,\D)$ is a
Hopf algebra is called CQT Hopf algebra.

\begin{proposition}[cf. Proposition $1.3.1.b$ of \cite{CQG}]
\label{prop1.7}
Let $\mbox{\em Poly}(G) = ({\cal B},\D)$ and $({\cal B},\D,{\cal R})$ be a
CQT Hopf algebra.  Then $(R^{vw})^{-1}$ exist for any $v,w \in \mbox{ {\em
Rep} }G$.  Moreover ${\cal R}' = {\cal R}(S \otimes id)$ satisfies
\[
{\cal R}'(x^{(1)} \otimes y^{(1)}){\cal R}(x^{(2)} \otimes y^{(2)}) = {\cal
R}(x^{(1)} \otimes y^{(1)}){\cal R}'(x^{(2)} \otimes y^{(2)}) = (\e \otimes
\e)(x \otimes y)
\]
(i.e. ${\cal R}' = {\cal R}^{-1}$) and
\[
{\cal
R}'(v_{il} \otimes w_{jk}) = (R^{vw})^{-1}_{ij,kl},\ i,l = 1,\dots,\dim v,\ j,k
= 1,\dots,\dim w.
\]
\end{proposition}

We say that $({\cal B},\D)$ is a ${}^*$-bialgebra if $({\cal B},\D)$ is a
bialgebra, ${\cal B}$ is a ${}^*$-algebra and
\begin{equation}
\label{eq1.16}
(* \otimes *)\D(x) = \D(x^*)
\end{equation}
for $x \in {\cal B}$.
A Hopf algebra which is a ${}^*$-bialgebra is called Hopf
${}^*$-algebra.

\begin{proposition}
\label{prop1.8}
Let $({\cal B},\D)$ be the bialgebra defined in Proposition~\ref{prop1.1}.
Suppose there exists an involution $\sim: {\cal J} \to {\cal J}$ such that
${\tilde 0} = 0$, $d_{\tilde\a}=d_{\a}$
and ${\tilde W}_m \in \mbox{ {\em Mor}}(w^{{\tilde
\a}_{ms_m}}
\otimes \dots \otimes w^{{\tilde \a}_{m1}},w^{{\tilde \b}_{mt_m}} \otimes
\dots
\otimes w^{{\tilde \b}_{m1}})$ where $({\tilde W}_m)_{i_{t_m}\dots i_1,j_{s_m}
\dots j_1} = \overline{(W_m)_{i_1\dots i_{t_m},j_1\dots j_{s_m}}}$.  Then there
exists a unique $*: {\cal B} \to {\cal B}$ such that $({\cal B},\D)$ is a
${}^*$-bialgebra and $\overline{w^{\a}} = w^{\tilde \a}$, $\a \in {\cal J}$
($-$
was defined in the Introduction).
\end{proposition}

\medskip
\noindent
{\bf Proof}.  Uniqueness is trivial.  Our assumptions imply that $z^{\a}_{ij} =
w^{\tilde \a}_{ij}$ satisfy (\ref{eq1.1})--(\ref{eq1.2}) in the conjugate
algebra ${\cal B}^j$ (conjugate vector space and opposite multiplication).
Therefore the desired $*$ exists (we check $*^2 = id$ and (\ref{eq1.16}) on the
generators $w^{\a}_{ij}$). \qed

\begin{proposition}[cf. the proofs of Proposition $1.3.d.ii$ and Theorem
$1.4.6$ of \cite{CQG}]\mbox{}\linebreak
\label{prop1.9}
Let $({\cal B},\D,{\cal R})$ be a CQT bialgebra and $\mbox{\em Poly}(G) =
({\cal B},\D)$ be a
${}^*$-bialgebra.
Suppose
$M \subset \mbox{ {\em Rep} } G$ is such that the
matrix elements of representations from $M$ generate ${\cal B}$ as
unital algebra.
  The following are equivalent:

{\em 1)}  $\overline{{\cal R}(y^* \otimes x^*)} = {\cal R}(x \otimes y)$, $x,y
\in {\cal B}$

{\em 2)}
\begin{equation}
\label{eq1.17}
\overline{(R^{{\bar w}{\bar v}})_{ji,lk}} = R_{ij,kl}^{vw},\ j,k =
1,\dots,\dim v,\ i,l = 1,\dots,\dim w,\ v,w \in M.
\end{equation}
If one of the above conditions is satisfied, $({\cal B},\D,{\cal
R})$ is called CQT ${}^*$-bialgebra.  If moreover $({\cal B},\D)$ is a Hopf
algebra,
$({\cal B},\D,{\cal R})$ is called CQT Hopf ${}^*$-algebra (see
Definition~$1.1$
of
\cite{CQG}).
\end{proposition}

\begin{proposition}[cf. \cite{Lyu}]
\label{prop1.10}
Let $\mbox{\em Poly}(G) = ({\cal B},\D)$ be a bialgebra and ${\cal R} \in
({\cal B}
\otimes {\cal B})'$.
Suppose
$M \subset \mbox{ {\em Rep} } G$ is such that the
matrix elements of representations from $M$ generate ${\cal B}$ as
unital algebra.
The following are equivalent:

{\em 1)}  ${\cal R}(x^{(1)} \otimes y^{(1)}){\cal R}(y^{(2)} \otimes x^{(2)}) =
(\e \otimes \e)(x \otimes y)$, $x,y \in {\cal B}$ (i.e. ${\cal R}_{21} = {\cal
R}^{-1}$)

{\em 2)}
\begin{equation}
\label{eq1.18}
(R^{vw})^{-1} = R^{wv},\ v,w \in M.
\end{equation}
If one of the above conditions is satisfied, we replace CQT by CT
(cotriangular) in all definitions (cf. \cite{Dr}).
\end{proposition}

\medskip
\noindent
{\bf Proof}.  Using (\ref{eq1.7})--(\ref{eq1.9}), one can assume that (in the
condition 2)) matrix elements of representations from $M$ linearly span
${\cal B}$.  Then in 1) it is enough to consider $x = v_{ij}$, $y = w_{kl}$,
$i,j = 1,\dots,\dim v$, $k,l = 1,\dots,\dim w$, $v,w \in M$. Then 1) is
equivalent to (\ref{eq1.18}). \qed

\section{Quasitriangular structures on quantum Lorentz groups}
\label{sec2}

In this section we classify coquasitriangular (${}^*$-) structures on quantum
Lorentz groups of \cite{WZ}.  For the quantum Lorentz group of
\cite{L} examples of such structures were given in \cite{Ta} and
Remark~7 of
Section~3 of
\cite{CQG}.  We also classify (cf. \cite{WP}) coquasitriangular (${}^*$-)
structures on (complex) $SL_q(2)$ groups and their real forms.

We first recall the definition of Hopf ${}^*$-algebra corresponding to a
quantum
Lorentz group \cite{WZ} essentially repeating the arguments of Theorem 1.1
of \cite{WZ}.
  The bialgebra structure of $({\cal A},\D)$ is obtained
by the construction of Proposition~\ref{prop1.1} with $\{w^{\a}: \a \in {\cal
J}\} = \{w^0,w,{\tilde w}\}$.
Here the relations (\ref{eq1.1}) mean that $w^0 =
(I_{\cal B})$ and the relations (\ref{eq1.2}) are given by
\setcounter{equation}{0}
\begin{equation}
\label{eq2.1}
(w \otimes w)E = Ew^0,
\end{equation}
\begin{equation}
\label{eq2.2}
({\tilde w} \otimes {\tilde w}){\tilde E} = {\tilde E}w^0,
\end{equation}
\begin{equation}
\label{eq2.3}
E'(w \otimes w) = w^0E',
\end{equation}
\begin{equation}
\label{eq2.4}
{\tilde E}'({\tilde w} \otimes {\tilde w}) = w^0{\tilde E}',
\end{equation}
\begin{equation}
\label{eq2.5}
X(w \otimes {\tilde w}) = ({\tilde w} \otimes w)X,
\end{equation}
where
\begin{equation}
\label{eq2.6}
{\tilde E} = \tau {\bar E},\ {\tilde E}' = {\bar E}'\tau,
\end{equation}
$E \in \mbox{ Lin}({\cb},{\cb}^2 \otimes {\cb}^2)$, $E' \in \mbox{
Lin}({\cb}^2
\otimes {\cb}^2,{\cb})$, $X \in \mbox{ Lin}({\cb}^2 \otimes {\cb
}^2,\
{\cb}^2 \otimes {\cb}^2)$ satisfy
\begin{equation}
\label{eq2.7}
(E' \otimes 1\!\!1)(1\!\!1 \otimes E) = 1\!\!1,
\end{equation}
\begin{equation}
\label{eq2.8}
(X \otimes 1\!\!1)(1\!\!1 \otimes X)(E \otimes 1\!\!1) = 1\!\!1 \otimes E,
\end{equation}
\begin{equation}
\label{eq2.9}
\tau {\bar X} \tau = \b^{-1}X,
\end{equation}
$E'E \ne 0$, $X$ is invertible, $\b \in {\cb}{\backslash}\{0\}$.

{\bf Warning}:  Our choice of $X$ may differ from the choice of \cite{WZ} by a
multiplicative nonzero constant.

Moreover, (\ref{eq2.7}) implies that $E$ and $E'$ are inverse one to another as
matrices,
\begin{equation}
\label{eq2.10}
(1\!\!1 \otimes E')(E \otimes 1\!\!1) = 1\!\!1.
\end{equation}
Hence, $E$ is left nondegenerate, $E'$ is right nondegenerate, $w^{-1}$ exists
(see Proposition~\ref{prop1.5}).  Similarly,
$(1\!\!1 \otimes {\tilde E}')({\tilde E} \otimes 1\!\!1) = 1\!\!1$,
${\tilde E}$
is left nondegenerate,
${\tilde E}'$ is right nondegenerate, ${\tilde w}^{-1}$ exists.  But
$(w^0)^{-1} = w^0$ and Proposition~\ref{prop1.6} implies that $({\cal A},\D)$
is a
Hopf algebra.

Setting $w^{\sim} = {\tilde w}$, ${\tilde w}^{\sim} = w$, $w^0{}^{\sim} = w^0$,
the assumptions of Proposition~\ref{prop1.8} are satisfied and $({\cal A},\D)$
becomes a ${}^*$-bialgebra where $*$ is defined by ${\bar w} = {\tilde w}$.  It
has the same Poincar\'e series as the classical $SL(2,{\cb})$ group
(Theorem~$1.2$ of \cite{WZ}).  We may assume that

\begin{enumerate}
\item $E = e_1 \otimes e_2 - qe_2 \otimes e_1$, $q \in {\cb
}{\backslash}\{0,i,-i\}$, $X = \a\tau Q$, $Q$ is given by (13)--(19) of
\cite{WZ} ($q=-1$ in (17)--(19)),
 $\a = t^{-1/2}$ for (13), $\a = (-t)^{-1/2}$ for (14), $\a =
q^{1/2}$ for (15), $\a = (-q)^{1/2}$ for (16), $\a = (s^2 - 1)^{-1/2}$ for
(17), $q = (p^2-1)^{-1/2}$ for (18), $\a = 1/2$ for (19), or
\item $E = e_1 \otimes e_2 - e_2 \otimes e_1 + e_1 \otimes e_1$ (in that
case we set $q = 1$),
$X = \tau Q$, $Q$ is given by (20)--(21) of \cite{WZ},
\end{enumerate}
$e_1,e_2$ is the canonical basis of ${\cb}^2$.  Moreover, $\b =
t/|t|$ for (13)--(14), $\b = q/|q|$ for (15), $\b = -q/|q|$
for (16),
$\b = -i \mbox{ sgn  Im}(s)$ for (17), $\b = \mbox{ sgn}(|p| - 1)$ for
(18), $\b = 1$ for (19)--(21).  In all cases $\b^4 = 1$.

We shall find all ${\cal R} \in ({\cal A} \otimes {\cal A})'$ such that $({\cal
A},\D,{\cal R})$ is a CQT Hopf algebra.  So (cf. Theorem~\ref{th1.4}) we need
to determine $R^{ww}$, $R^{{\bar w}{\bar w}}$, $R^{{\bar w}w}$, $R^{w{\bar w}}$
such that 20 relations (\ref{eq1.12})--(\ref{eq1.13}) and 4 relations
(\ref{eq1.14}) are satisfied (we assume (\ref{eq1.11})).  We shall use them in
the following.  Irreducibility of $w \otimes {\bar w}$ (see
\cite{WZ}) and (\ref{eq2.5}) give $R^{w{\bar w}} = \e_XX$, $R^{{\bar w}w} =
\e'_XX^{-1}$, where $\e_X,\e'_X \in {\cb}{\backslash}\{0\}$ (cf.
Proposition~\ref{prop1.7}).  Moreover, $D = R^{ww}$ and ${\tilde D} = \b
R^{{\bar w}{\bar w}}$ must satisfy $(2.11)$ and
$(2.17)$--$(2.20)$ of \cite{POI} (with $L,{\tilde L}$ replaced by $D,{\tilde
D}$), hence they are given by $(2.21)$--$(2.22)$ of \cite{POI}, i.e. $R^{ww} =
L_i$, $R^{{\bar w}{\bar w}} = \tau\overline{L_j^{-1}}\tau$, $i,j = 1,2,3,4$,
$L_i = q_i(1\!\!1 + q_i^{-2}EE')$, $q_{1,2} = \pm q^{1/2}$, $q_{3,4} = \pm
q^{-1/2}$.  Using $(2.1)$, $(2.3)$--$(2.4)$ of \cite{POI}, we obtain $\e_X =
\pm 1$, $\e'_X = \pm 1$, $\b = \pm 1$.  After some calculations one gets that
the 24 relations are satisfied.  So we get $4\cdot
4\cdot 4$
${\cal R}$ for $q \ne \pm 1$, $\b = \pm 1$, $2\cdot 2\cdot 4$ ${\cal R}$ for $q
= \pm 1$, $\b = \pm 1$ and no ${\cal R}$ for $\b = \pm i$.

Set $M = \{w,{\bar w}\}$.  According to Proposition~\ref{prop1.9}, we get a CQT
Hopf ${}^*$-algebra iff 4 relations (\ref{eq1.17}) are satisfied iff $\b = 1$
and
$q_j = q_i^{-1}$ ($4\cdot 4$ ${\cal R}$ for $q \ne \pm 1$, $\b = 1$, $2\cdot 4$
${\cal R}$ for $q = \pm 1$, $\b = 1$, no ${\cal R}$ for $\b \ne 1$).  According
to Proposition~\ref{prop1.10}, we get CT Hopf algebra iff 4 relations
(\ref{eq1.18}) are satisfied iff $q = 1$,
$\e'_X=\e_X$ ($2\cdot 2\cdot 2$ ${\cal R}$ for $q =
1$, $\b = \pm 1$, no ${\cal R}$ otherwise).  Clearly, we get CT Hopf
${}^*$-algebra iff $q = \b = 1$ and $q_j = q_i^{-1}$, $\e'_X=\e_X$
 ($2\cdot 2$ ${\cal R}$ for $q
= \b = 1$, no ${\cal R}$ otherwise).

Let us also recall $SL_q(2)$ groups \cite{W1}, \cite{Dr}.  The bialgebra
structure is obtained by the construction of Proposition~\ref{prop1.1} with
$\{w^{\a}: \a \in {\cal J}\} = \{w^0,w\}$, the relations (\ref{eq1.1}) mean
that $w^0 = (I_{\cal B})$ and those of (\ref{eq1.2}) are
given by (\ref{eq2.1}),
(\ref{eq2.3}), where $E,E'$ are as above with $q \in {\cb}{\backslash}\{0\}$
(only the case~1.).  One gets that $({\cal A},\D)$ is a Hopf algebra.

We shall describe (cf. \cite{WP})
all ${\cal R} \in ({\cal A} \otimes {\cal A})'$
such
that
$({\cal A},\D,{\cal R})$ is a CQT Hopf algebra.  Due to Theorem~\ref{th1.4} we
should find $D = R^{ww} \in \mbox{ Mor}(w \otimes w,w \otimes w)$ satisfying 4
relations (\ref{eq1.12})--(\ref{eq1.13}).  This means $D = L_i$, $i =
1,2,3,4$, so we get 4 ${\cal R}$ for $q \ne \pm 1$, $0$, and $2$ ${\cal R}$ for
$q = \pm 1$.  In other words $R^{ww}
 = \pm L_1^{\pm 1}$ where $L_1 = q^{1/2}(1\!\!1
+ q^{-1}EE')$.  We get CT Hopf algebras iff $L_1^2 = 1\!\!1$ (see Proposition
$\ref{prop1.10}$ with $M = \{w\}$) iff $q = 1$.

Let us pass to real forms.  Then $({\cal A},\D)$ becomes a Hopf ${}^*$-algebra
where
$*$ is defined by ${\bar w} = v^c$, $v = BwB^{-1}$, $B = 1\!\!1$ for $SU_q(2)$,
$B = \mbox{ diag}(1,-1)$ for $SU_q(1,1)$, $q \in {\rb}{\backslash}\{0\}$
(cf. \cite{W1}, \cite{RTF}, \cite{CQG}), ${\bar w} = w$ for $SL_q(2,{\rb})$,
$|q| = 1$ (cf. \cite{RTF} and Proposition~\ref{prop1.8}).  For $SU_q(2)$,
$SU_q(1,1)$ we get CQT Hopf ${}^*$-algebras iff $L_1$ is hermitian w.r.t. the
canonical scalar product in ${\cb}^2$ (see Proposition $\ref{prop1.9}.2$
with $M =
\{w\}$ and the proof of Theorem $1.4.6$ of \cite{CQG}) iff $q > 0$.  For
$SL_q(2,{\rb})$ we get CQT Hopf ${}^*$-algebras iff $\tau {\bar L}_1\tau = L_1$
(see Proposition $\ref{prop1.9}.2$
with $M = \{w\}$) iff $q = 1$.  So we have CT Hopf
${}^*$-algebras iff $q = 1$ (for each real form).

\section{Quasitriangular structures on inhomogeneous quantum groups}
\label{sec3}

For any Hopf algebra $\mbox{Poly}(H) = ({\cal A},\D)$ satisfying certain
properties one can construct \cite{INH} a Hopf algebra $\mbox{Poly}(G) = ({\cal
B},\D)$ which describes an inhomogeneous quantum group.  For certain CQT Hopf
algebra structures $({\cal A},\D,{\cal R}_{\cal A})$ we find all CQT Hopf
algebra
structures $({\cal B},\D,{\cal R})$ such that ${\cal R}_{{|}_{{\cal
A} \otimes {\cal A}}} = {\cal R}_{\cal A}$.  The ${}^*$-structures and
cotriangularity are studied as well.  In particular we find all C(Q)T Hopf
(${}^*$-) algebra structures on quantum Poincar\'e groups.

Throughout the Section $\mbox{Poly}(H) = ({\cal A},\D)$ is any bialgebra such
that
\begin{enumerate}
\begin{enumerate}
\item  each representation of $H$ is completely reducible,
\item $\Lambda$ is an irreducible representation of $H$,
\item $\mbox{Mor}(v \otimes w,\Lambda \otimes v \otimes w) = \{0\}$ for any two
irreducible representations $v,w$ of $H$.
\end{enumerate}
\end{enumerate}

Moreover, we assume that $f_{ij}$, $\eta_i \in {\cal A}'$, $i,j = 1,\dots$, $N
= \dim \Lambda$, are given and satisfy
\begin{enumerate}
\item ${\cal A} \ni a \to \rho(a) = \left( \begin{array}{cc} f(a) & \eta(a) \\
0 & \e(a) \end{array} \right) \in M_{N+1}({\cb})$ is a unital homomorphism,
\item $\Lambda_{st}(f_{tr} * a) = (a * f_{st})\Lambda_{tr}$ for $a \in
{\cal A}$,
\item $R^2 = 1\!\!1$ where $R_{ij,sm} = f_{im}(\Lambda_{js})$,
\item $(\Lambda \otimes \Lambda)_{kl,ij}(\tau^{ij} * a) = a * \tau^{kl}$ for $a
\in {\cal A}$ where $\tau^{ij} = (R - 1)_{ij,mn}(\eta_n * \eta_m -
\eta_m(\Lambda_{ns})\eta_s + T_{mn} \e -f_{nb} * f_{ma} T_{ab})$,
\item $A_3{\tilde F} = 0$ where $A_3 = 1\!\!1 \otimes 1\!\!1 \otimes 1\!\!1 - R
\otimes 1\!\!1 - 1\!\!1 \otimes R + (R \otimes 1\!\!1)(1\!\!1 \otimes R) +
(1\!\!1 \otimes R)(R \otimes 1\!\!1) - (R \otimes 1\!\!1)(1\!\!1 \otimes R)(R
\otimes 1\!\!1)$, ${\tilde F}_{ijk,m} = \tau^{ij}(\Lambda_{km})$,
\item $A_3(Z \otimes 1\!\!1 - 1\!\!1 \otimes Z)T = 0$,
$RT=-T$ where $Z_{ij,k} =
\eta_i(\Lambda_{jk})$.
\end{enumerate}
In particular, 4. and 5. are satisfied if $\tau^{ij} = 0$.

The inhomogeneous quantum group $G$ corresponds to the
bialgebra $\mbox{Poly}(G) =
({\cal B},\D)$ defined  (cf. Corollary $3.8.a$ of \cite{INH}) as follows:
${\cal B}$ is the universal unital algebra generated by ${\cal A}$ and $y_i$,
$i = 1,\dots,N$, with the relations $I_{\cal B} = I_{\cal A}$,
\setcounter{equation}{0}
\begin{equation}
\label{eq3.0}
y_sa = (a * f_{st})y_t + a * \eta_s - \Lambda_{st}(\eta_t * a),\ a \in {\cal
A},
\end{equation}
\begin{equation}
\label{eq3.0'}
(R - 1\!\!1)_{kl,ij}(y_iy_j - \eta_i(\Lambda_{js})y_s + T_{ij} -
\Lambda_{im}\Lambda_{jn}T_{mn}) = 0.
\end{equation}
Moreover, $({\cal A},\D)$ is a subbialgebra of $({\cal B},\D)$ and $\D y_i =
\Lambda_{ij} \otimes y_j + y_i \otimes I$ ($y_i$ were denoted by $p_i$ in
\cite{INH}).  In particular,
 ${\cal P} = \left( \begin{array}{cc} \Lambda & y \\
0 & I \end{array} \right)$ is a representation of $G$.

\begin{remark}
\label{rem3.1}
{\em If $H$ is a matrix group, $\Lambda$ its fundamental representation and
$({\cal A},\D)$ its corresponding Hopf algebra (generated by $\Lambda_{ij}$,
$\Lambda^{-1}_{ij}$) then, assuming (a)--(c) and setting $f_{ij} = \d_{ij}\e$,
$\eta_i = 0$, $T = 0$, $({\cal B},\D)$ corresponds to
\[
G = H \sd {\rb}^N
= \left\{ g=\left( \begin{array}{cc} h & a \\ 0 & 1 \end{array} \right) \in
M_{N+1}({\cb}): h \in H,\ a \in {\rb}^N\right\},
\]
$f(g)=f(h)$, $y_i(g)=a_i$, $f\in {\cal A}$, $i=1,\ldots,N$, $g\in G$.}
\end{remark}

According to Corollary $3.8.b$ and the proof of Proposition $3.12$ of
\cite{INH}, the bialgebra $({\cal B},\D)$ can be obtained by the construction
of
Proposition~\ref{prop1.1} with $\{w^{\a}: \a \in {\cal J}\} = \mbox{ Irr } H
\cup \{{\cal P}\}$.
Here the relations (\ref{eq1.1}) mean that $w^0 \equiv (I_{\cal
A}) = (I_{\cal B})$ and the relations (\ref{eq1.2}) are given by
\begin{equation}
\label{eq3.1}
({\cal P} \otimes {\cal P})R_P = R_P({\cal P} \otimes {\cal P}),
\end{equation}
\begin{equation}
\label{eq3.2}
({\cal P} \otimes w)N_w = N_w(w \otimes {\cal P}),\ w \in \mbox{ Irr } H,
\end{equation}
\begin{equation}
\label{eq3.3}
(w \otimes w')S_{ww'w''}^{\a} = S_{ww'w''}^{\a}w'',\ w,w',w'' \in \mbox{ Irr }
H,\ \a = 1,\dots,c_{ww'}^{w''},
\end{equation}
\begin{equation}
\label{eq3.4}
{\cal P}i = i\Lambda,
\end{equation}
\begin{equation}
\label{eq3.5}
s{\cal P} = w^0s,
\end{equation}
where
\begin{equation}
\label{eq3.5'}
R_P = \left( \begin{array}{cccc} R & Z & -R \cdot Z & (R - 1\!\!1)T \\
0 & 0 & 1\!\!1 & 0 \\
0 & 1\!\!1 & 0 & 0 \\
0 & 0 & 0 & 1
\end{array} \right),\ N_w = \left( \begin{array}{c} G_w,\ H_w \\
0,\ 1\!\!1
\end{array} \right),
\end{equation}
$(G_w)_{iC,Dj} = f_{ij}(w_{CD})$, $(H_w)_{iC,D} = \eta_i(w_{CD})$, $w\in\mbox{
 Rep } H$,
$R =
G_{\Lambda}$, $Z = H_{\Lambda}$, $S_{ww'w''}^{\a}$ $(\a = 1,
\dots,c_{ww'}^{w''})$
form a basis of $\mbox{Mor}(w'',w \otimes w')$, $i: {\cb}^N \to {\cb}^N
\oplus {\cb}$, $s: {\cb}^N \oplus {\cb} \to {\cb}$ are the
canonical mappings.  In the following we assume that $({\cal A},\D)$ is a Hopf
algebra.  Then (Proposition $3.12$ of \cite{INH}) $({\cal B},\D)$ is also a
Hopf algebra and $G_w^{-1}$ exist:
\begin{equation}
\label{eq3.6}
(G_w^{-1})_{Ak,lB} = f_{kl}(w_{AB}^{-1}).
\end{equation}
If $({\cal A},\D)$ is a Hopf ${}^*$-algebra then we also assume ${\bar
\Lambda} =
\Lambda$,
\begin{equation}
\label{eq3.6a}
f_{ij}(S(a^*)) = \overline{f_{ij}(a)},\ \eta_i(S(a^*)) =
\overline{\eta_i(a)},\
i,j = 1,\dots,N,\ a \in {\cal A},
\end{equation}
${\tilde T} - T \in \mbox{ Mor}(w^0,\Lambda \otimes \Lambda)$, where ${\tilde
T}_{ij} = \overline{T_{ji}}$.  Then \cite{INH} $({\cal B},\D)$ has a unique
Hopf ${}^*$-algebra structure such that $({\cal A},\D)$ is its Hopf
${}^*$-subalgebra
and $y_i^* = y_i$.

In the following we assume $\mbox{ Mor }(I,\Lambda\otimes\Lambda)\cap
\mbox{ ker }(R+1\!\!1)=\{0\}$, i.e. $\mbox{ Mor }(I,\Lambda\otimes\Lambda)
\subset\mbox{ ker }(R-1\!\!1)$. Then (using (4.14) of \cite{INH}) $\tilde T=T$.
The main result of the present paper is contained in

\begin{theorem}
\label{th3.2}
Let $\mbox{\em Poly}(H) = ({\cal A},\D)$, $\mbox{\em Poly}(G) = ({\cal B},\D)$
be as above, $({\cal A},\D,{\cal R}_{\cal A})$ be a CQT Hopf algebra such that
\begin{equation}
\label{eq3.7}
R^{v\Lambda} = c_vG_v,\ R^{\Lambda v} = c'_v G_v^{-1},\ v \in
\mbox{ {\em Irr} } H,
\end{equation}
$c_v,c'_v \in {\cb}{\backslash}\{0\}$.  We are interested in CQT Hopf
algebra structures $({\cal B},\D,{\cal R})$ such that
\begin{equation}
\label{eq3.7'}
{\cal R}_{{|}_{{\cal A} \otimes {\cal A}}} = {\cal R}_{\cal A}.
\end{equation}
One has:
\begin{enumerate}
\item Such a structure exists iff
\begin{equation}
\label{eq3.8}
\tau^{ij} = 0,\ i,j = 1,\dots,N,
\end{equation}
\begin{equation}
\label{eq3.9}
R^{v\Lambda} = G_v,\ R^{\Lambda v} = G_v^{-1},\ v \in \mbox{ {\em Irr } } H.
\end{equation}
\item Suppose (\ref{eq3.8})--(\ref{eq3.9}).  Then such structures are in one to
one correspondence with $m \in \mbox{ {\em Mor}}(w^0,\Lambda \otimes \Lambda)$
satisfying
\begin{equation}
\label{eq3.9''}
(f_{jb}*f_{ia})m_{ab}=m_{ij}\e,\quad i,j=1,\ldots,N,
\end{equation}
and are determined by
\begin{equation}
\label{eq3.9'}
R^{vw} = R^{vw} \mbox{ for } {\cal R}_{\cal A},\ v,w \in \mbox{ {\em Irr} } H,
\end{equation}
\begin{equation}
\label{eq3.10}
R^{v{\cal P}} = N_v,\ R^{{\cal P}v} = N_v^{-1},\ v \in \mbox{ {\em Irr} } H,
\end{equation}
\begin{equation}
\label{eq3.11}
R^{{\cal P}{\cal P}} = R_P + m_P
\end{equation}
where
\begin{equation}
\label{eq3.12}
m_P = \left( \begin{array}{cccc} 0 & 0 & 0 & m \\
0 & 0 & 0 & 0 \\
0 & 0 & 0 & 0 \\
0 & 0 & 0 & 0
\end{array} \right) .
\end{equation}
\item Let ${\cal R}$ be as in $2.$
and let $({\cal A},\D)$, $({\cal B},\D)$ be Hopf ${}^*$--algebras as in
the text before the Theorem.  Then $({\cal B},\D,{\cal R})$ is a CQT
Hopf ${}^*$-algebra iff $({\cal A},\D,{\cal R}_{\cal A})$ is a CQT Hopf
${}^*$-algebra and
\begin{equation}
\label{eq3.12a}
m_{ij}= \overline{m_{ji}},\quad i,j=1,\ldots,N.
\end{equation}
\item Let ${\cal R}$ be as in $2$.  Then $({\cal B},\D,{\cal R})$ is a CT Hopf
algebra iff $({\cal A},\D,{\cal R}_{\cal A})$ is a CT Hopf algebra and $m =
0$.
\end{enumerate}
\end{theorem}

\medskip
\noindent
{\bf Proof}.  {\bf Ad 1--2}.  Each such structure is (see Theorem~\ref{th1.4})
uniquely determined by $R^{vw}$, $R^{v{\cal P}}$, $R^{{\cal P}v}$ and $R^{{\cal
P}{\cal P}}$ satisfying (\ref{eq1.12})--(\ref{eq1.14}) (we assume
(\ref{eq1.11})), $v,w \in \mbox{ Irr } H$.  Using (\ref{eq3.7'}), we get
(\ref{eq3.9'}).  In virtue of the properties of ${\cal R}_{\cal A}$ the formula
(\ref{eq1.14}) for $R^{vw}$ follows.  Moreover, (\ref{eq1.12})--(\ref{eq1.13})
for
(\ref{eq3.4}) and $w^{\g} = v \in \mbox{ Irr } H$ mean $R^{v{\cal P}}(1\!\!1
\otimes i) = (i \otimes 1\!\!1)R^{v\Lambda}$, $R^{{\cal P} v}(i \otimes
1\!\!1) =
(1\!\!1 \otimes i)R^{\Lambda v}$.  That and (\ref{eq3.7}) give
\[
R^{v{\cal P}} = c_v \left( \begin{array}{c}
G_v,\ ? \\
0,\ ?
\end{array} \right) = c_vN_v,\ R^{{\cal P} v} = c'_v \left( \begin{array}{cc}
G_v^{-1} & ? \\
0 & ?
\end{array} \right) = c'_vN_v^{-1}
\]
where the second equalities follow from (\ref{eq1.14}) for $R^{v{\cal
P}},R^{{\cal P}v}$, $(\ref{eq3.2})$, the independence of $1,y_i$ $(i =
1,\dots,N)$ over ${\cal A}$ (in left and also in right module, see
Corollary 3.6 and (1.4) of \cite{INH}) and the
condition (c).  Using (\ref{eq1.12})--(\ref{eq1.13}) for (\ref{eq3.5}) and
$w^{\g} = v \in \mbox{ Irr } H$, one obtains $c_v = c'_v = 1$, we get
(\ref{eq3.9}), (\ref{eq3.10}).  In virtue of (\ref{eq1.12})--(\ref{eq1.13}) for
(\ref{eq3.4})--(\ref{eq3.5}) and $w^{\g} = {\cal P}$
\[
R^{{\cal P}{\cal P}} = \left( \begin{array}{cccc}
R & Z & -R\cdot Z & ? \\
0 & 0 & 1\!\!1 & 0 \\
0 & 1\!\!1 & 0 & 0 \\
0 & 0 & 0 & 1
\end{array} \right) = R_P + m_P
\]
for some $m \in \mbox{ Mor}(w^0,\Lambda \otimes \Lambda)$
where the second equality uses (\ref{eq1.14}) for $R^{{\cal P}{\cal P}}$,
$(\ref{eq3.1})$ and (\ref{eq3.12}).  Thus (\ref{eq3.11}) follows.  We
set $R_Q = R_P + m_P = R^{{\cal P}{\cal P}}$ and replace (\ref{eq3.1}) by
equivalent (assuming (\ref{eq3.3})) relation
\begin{equation}
\label{eq3.13}
({\cal P} \otimes {\cal P})R_Q = R_Q({\cal P} \otimes {\cal P}).
\end{equation}
The relations (\ref{eq1.12})--(\ref{eq1.13}) for
(\ref{eq3.13}),
(\ref{eq3.2})
 and any $w^{\g} \in \mbox{ Irr } H \cup \{{\cal P}\}$ are
equivalent to
\begin{equation}
\label{eq3.14}
(R_Q \otimes 1\!\!1)(1\!\!1 \otimes R_Q)(R_Q \otimes 1\!\!1) = (1\!\!1 \otimes
R_Q)(R_Q \otimes 1\!\!1)(1\!\!1 \otimes R_Q),
\end{equation}
\begin{equation}
\label{eq3.15}
(1\!\!1 \otimes N_v)(N_v \otimes 1\!\!1)(1\!\!1 \otimes R_Q) = (R_Q \otimes
1\!\!1)(1\!\!1 \otimes N_v)(N_v \otimes 1\!\!1),\ v \in \mbox{ Irr } H,
\end{equation}
\begin{equation}
\label{eq3.16}
(1\!\!1 \otimes R^{vw})(N_v \otimes 1\!\!1)(1\!\!1 \otimes N_w) = (N_w \otimes
1\!\!1)(1\!\!1 \otimes N_v)(R^{vw} \otimes 1\!\!1),\ v,w \in \mbox{ Irr } H.
\end{equation}
According to Proposition $3.14$ of \cite{INH} and its proof, (\ref{eq3.14}) is
equivalent to ${\tilde F} = 0$.  Let us denote the standard basis elements in
${\cb}^{\dim v}$, $v \in \mbox{ Irr } H$, by $h_i^v$, $i = 1,\dots,\dim v$,
$e_i=h^{\Lambda}_i$ and in $\cb$ by $f=1$.
Using $(3.65)$ of \cite{INH} and
\[
\begin{array}{rll}
N_v(h_i^v \otimes e_j) &= &(G_v)_{kl,ij} e_k \otimes h_l^v, \\
N_v(h_i^v \otimes f) &= &f \otimes h_i^v + (H_v)_{kl,i} e_k \otimes h_l^v,
\end{array}
\]
we find that (\ref{eq3.15}) on $h_i^v \otimes e_j \otimes e_k$ follows from the
last formula before Proposition $3.14$ in \cite{INH}, on $h_i^v \otimes e_j
\otimes f$, $h_i^v \otimes f \otimes e_j$ follows from $(2.18)$ of \cite{INH},
on $h_i^v \otimes f \otimes f$ is equivalent
(using $Rm=m$) to $\tau^{ij}(v_{AB}) = 0$
and (\ref{eq3.9''}) applied to $v_{AB}$
 for all
$i,j,A,B$.  So (\ref{eq3.15}) is equivalent to (\ref{eq3.8}) (which implies
${\tilde F} = 0$)
and (\ref{eq3.9''}).
  We also get that (\ref{eq3.16}) on $h_i^v \otimes h_j^w
\otimes e_k$ follows from $(1\!\!1 \otimes R^{vw})(R^{v\Lambda} \otimes
1\!\!1)(1\!\!1 \otimes R^{w\Lambda}) = (R^{w\Lambda} \otimes 1\!\!1)(1\!\!1
\otimes R^{v\Lambda})(R^{vw} \otimes 1\!\!1)$ (it can be obtained using
(\ref{eq1.10}), (\ref{eq1.6a}) for ${\cal R}_{\cal A} -$ cf.
\cite{Dr}, \cite{LT},
$(1.14)$ of \cite{CQG}), on $h_i^v \otimes h_j^w \otimes f$ follows
{}from the equality obtained by acting $\eta_i$ on (\ref{eq1.14}) for $R^{vw}$
and using
condition 1. (see the beginning of the Section).

The relations (\ref{eq1.12})--(\ref{eq1.13}) for (\ref{eq3.3}) and $w^{\g} \in
\mbox{ Irr } H$ follow from (\ref{eq1.6a}),
(\ref{eq1.8}) for ${\cal R}_{\cal A}$ while for $w^{\g} = {\cal
P}$ are equivalent to
\[
(N_w \otimes 1\!\!1)(1\!\!1 \otimes N_{w'})(S^{\a}_{ww'w''} \otimes 1\!\!1) =
(1\!\!1 \otimes S^{\a}_{ww'w''})N_{w''}.
\]
That on $h_i^{w''} \otimes e_j$ follows from (\ref{eq1.6a}), (\ref{eq1.8}) for
${\cal R}_{\cal A}$, on $h_i^{w''} \otimes f$ follows from the equality
obtained by acting $\eta_i$ on (\ref{eq3.3}).

\medskip
{\bf Ad 3}.  We need to check (\ref{eq1.17}) for $M = \mbox{ Irr } H
\cup
\{{\cal P}\}$.  For $R^{vw}$, $v,w \in \mbox{ Irr } H$, it is equivalent to the
fact that $({\cal A},\D,{\cal R}_{\cal A})$ is a CQT Hopf ${}^*$-algebra, for
$R^{{\cal P}v}$ and $R^{v{\cal P}}$ follows from (\ref{eq3.6}), (\ref{eq3.6a})
and the properties of $\eta_i$, for $R^{{\cal P}{\cal P}}$ it is equivalent
(using $(4.14)$ and the next formula of
\cite{INH}, $Rm=m$, $RT=-T$) to (\ref{eq3.12a}).

\medskip
{\bf Ad 4}. We need to check (\ref{eq1.18}) for $M = \mbox{ Irr } H \cup
\{{\cal P}\}$.  For $R^{vw}$, $v,w \in \mbox{ Irr } H$, it is equivalent to the
fact that $({\cal A},\D,{\cal R}_{\cal A})$ is a CT Hopf algebra, for $R^{{\cal
P}v}$ and $R^{v{\cal P}}$ follows from (\ref{eq3.10}), for $R^{{\cal P}{\cal
P}}$ it is equivalent (using $R_P^2 = 1\!\!1$) to $m = 0$. \qed

\medskip
{\bf Remark.} If the first formula of the condition 6. is replaced by
$0\neq A_3(Z\otimes 1\!\!1-1\!\!1\otimes Z)T\in\mbox{ Mor }
(I,\Lambda\otimes\Lambda\otimes\Lambda)$
(this is allowed by \cite{INH}) then (\ref{eq3.14})
is not satisfied and there is no CQT Hopf algebra structure on $({\cal B},\D)$.

\medskip
As an application we shall consider $({\cal A},\D) = \mbox{ Poly}(L)$ where $L$
is a quantum Lorentz group.  The corresponding inhomogeneous quantum groups are
called quantum Poincar\'e groups and are (almost) classified in \cite{POI}.
The classification of C(Q)T Hopf (${}^*$-) algebra structures on them
is given in

\begin{theorem}
\label{th3.3}
Let $\mbox{\em Poly}(P) = ({\cal B},\D)$ be the Hopf ${}^*$-algebra
corresponding
to a quantum Poincar\'e group $P$ \cite{POI} described by an admissible choice
of quantum Lorentz group (cases 1)--7)), $s = \pm 1$, $H$ and $T$.

\medskip
\noindent
{\em 1}.  Let us consider CQT Hopf algebra structures $({\cal B},\D,{\cal R})$
on $P$.  One has:

\medskip
{\em a)}  In the cases $1)$ (except $s=1$, $t=1$, $t_0 \ne 0$---
see Remark $1.8$ of
\cite{POI}\ {}\footnote{
In the old version of Remark 1.8 of \cite{POI} one should replace
$t$ by $t_0$ (except of expressions $t=1$). That $t_0$ is identified with $t$
of (3)--(4) of Ref. 16 of \cite{POI}.}
),
$2)$, $3)$, $4)$ (except $s=1$, $b \ne 0$), 5) (except $s=\pm1$,
$t=1$, $t_0\neq0$), 6), 7)
each such structure is uniquely determined by
\[
R^{ww} = kL,\ R^{w{\bar w}} = kX,\ R^{{\bar w}w} = qkX^{-1},\ R^{{\bar w}{\bar
w}} = qk{\tilde L}
\]
and (\ref{eq3.10})--(\ref{eq3.12}), where
\[
m = cm_0,\ m_0 = (V^{-1} \otimes V^{-1})(1\!\!1 \otimes X \otimes 1\!\!1)(E
\otimes \tau E), \ L = sq^{1/2}(1\!\!1+q^{-1}EE'),\ \tilde L=q\tau L\tau,
\]
where $s=\pm1$, $E,E'$ are fixed for fixed $P$ and given in \cite{POI},
$\ k = \pm 1$ (two possible ${\cal R}$ for each $c \in {\cb}$).

\medskip
{\em b)}  In the other cases there is no such structure.

\medskip
\noindent
{\em 2}.  Let ${\cal R}$ be as in $1$.  We get CQT Hopf ${}^*$-algebra iff
$q=1$ (which excludes the cases 5), 6), 7)) and $c \in
{\rb}$.

\medskip
\noindent
{\em 3}.  Let ${\cal R}$ be as in $1$.  We get CT Hopf algebra iff
$q=1$ (which excludes the cases 5), 6), 7)) and $c = 0$
(then it is also CT Hopf ${}^*$-algebra).
\end{theorem}

\medskip
\noindent
{\bf Proof}.  {\bf Ad 1}.  We shall use Theorem~\ref{th3.2}, the results of
\cite{POI} and Section~\ref{sec2}.  Thus $H$ is a quantum Lorentz group,
$\Lambda = V^{-1}(w \otimes {\bar w})V$ with $V_{CD,i} = (\s_i)_{CD}$ $(\s_0 =
1\!\!1,\s_1,\s_2,\s_3$ are the Pauli matrices), $q = \b = \pm 1$,
the assumptions about $H$ and $G$  are
satisfied.  Moreover, $G_w = (V^{-1} \otimes 1\!\!1)(1\!\!1 \otimes X)(L
\otimes 1\!\!1)(1\!\!1 \otimes V)$, $G_{\bar w} = (V^{-1} \otimes
1\!\!1)(1\!\!1 \otimes {\tilde L})(X^{-1} \otimes 1\!\!1)(1\!\!1 \otimes V)$,
where $L = sq^{1/2}(1\!\!1+q^{-1}EE')$
and $\tilde L=q\tau\overline{L^{-1}}\tau=
q\tau L\tau$.
The possible $\eta$ and $T$ are described in \cite{POI}.
According to the results of Section~\ref{sec2}, each CQT Hopf algebra structure
on $({\cal A},\D)$ is uniquely characterized by $R^{ww} = \e_LL$, $R^{{\bar
w}{\bar w}} = \e'_L{\tilde L}$, $R^{w{\bar w}} =
\e_XX$,
$R^{{\bar w}w}= \e'_XX^{-1}$,
 where $\e_L^2 = \e_L^{\prime 2} = \e_X^2 = \e_X^{\prime 2}=1$ (16 possible
${\cal R}_{\cal A}$).
Using (\ref{eq1.6a})--(\ref{eq1.6b}), (\ref{eq1.8})--(\ref{eq1.9}), one gets
(\ref{eq3.7}) for $v = w,{\bar w}$ with $c_w = \e_L\e_X$, $c_{\bar w} =
\e'_L\e'_X$, $c'_w = q\e_L\e'_X$, $c'_{\bar w} = q\e_X\e'_L$.
In virtue of
Proposition $2.1$ of \cite{POI} and
(\ref{eq1.6a})--(\ref{eq1.6b}), (\ref{eq1.8})--(\ref{eq1.9}) the condition
(\ref{eq3.7}) is satisfied for all $v \in \mbox{ Irr } H$.

Due to Proposition $3.13.2$ and Proposition $4.8$ of \cite{INH}, (\ref{eq3.8})
is equivalent to $\tau^{ij}(w_{AB}) = 0$, $i,j = 0,1,2,3$, $A,B = 1,2$, which
means (cf. the proof of Theorem $1.6$ of \cite{POI}) $\lambda = 0$, which
excludes the case 4), $s = 1$, $b \ne 0$, the case 1), $s = 1$,
$t=1$, $t_0 \ne 0$ and
the case 5), $s = \pm 1$, $t=1$, $t_0 \ne 0$ where $t_0\in\rb$ is
introduced in Remark $1.8$ of \cite{POI}.  Moreover, (\ref{eq3.9}) means that
$\e'_X=\e'_L=qk$, $\e_X=\e_L=k$ for some $k=\pm1$.
Using Theorem $\ref{th3.2}.1$--$2$,
$\mbox{Mor}(w^0,\Lambda \otimes \Lambda) = {\cb
}m_0$ and (\ref{eq3.9''}) for $m=m_0$ (it is enough to prove it on
$w_{AB}$, $w_{AB}{}^*$ when it follows from the 20 relations considered in
Section~{\ref{sec2}}), we get $1$.

{\bf Ad 2}.  We use $q=\b=1$ (which implies $q_j=q_i^{-1}$),
 $(m_0)_{ij} =
\overline{(m_0)_{ji}}$ and Theorem $\ref{th3.2}.3$.

{\bf Ad 3}.  We use $q = 1$ (which implies $\e'_X=\e_X$) and Theorem
$\ref{th3.2}.4$. \qed

\section{Enveloping algebras}
\label{sec4}

In this section we study enveloping algebras of inhomogeneous quantum groups.
We assume that $({\cal A},\D,{\cal R}_{\cal A})$ and $({\cal B},\D,{\cal R})$
are CQT Hopf algebras as in Theorem $\ref{th3.2}.1$--$2$ (e.g. as in Theorem
$\ref{th3.3}.1$).

We essentially
follow the scheme of \cite{RTF}
and \cite{UEA} but now we don't assume $Z=T=0$.
We define $l_{jl} \in {\cal B}'$, $j,l = 1,\dots,N,+$, by
\setcounter{equation}{0}
\begin{equation}
\label{eq4.0}
l_{jl}(x) = {\cal
R}(x \otimes {\cal P}_{jl})
\end{equation}
(in CT case $l$ corresponds to $L^{\pm}$ of \cite{RTF}
on the
subalgebra generated by ${\cal P}_{ac}$).  According to (\ref{eq1.5}) and
(\ref{eq1.10}) for $v = w = {\cal P}$,
\[
\begin{array}{rll}
R^{{\cal P}{\cal P}}_{ab,cd}l_{df}(x^{(1)})l_{ce}(x^{(2)}) &= &R^{{\cal P}{\cal
P}}_{ab,cd} {\cal R}(x^{(1)} \otimes {\cal P}_{df}){\cal R}(x^{(2)} \otimes
{\cal P}_{ce}) \\
&= &{\cal R}(x \otimes R^{{\cal P}{\cal P}}_{ab,cd}{\cal P}_{ce}{\cal P}_{df})
\\
&= &{\cal R}(x \otimes {\cal P}_{ac}{\cal P}_{bd}R^{{\cal P}{\cal P}}_{cd,ef})
\\
&= &{\cal R}(x^{(1)} \otimes {\cal P}_{bd}){\cal R}(x^{(2)} \otimes {\cal
P}_{ac})R^{{\cal P}{\cal P}}_{cd,ef} \\
&= &l_{bd}(x^{(1)})l_{ac}(x^{(2)})R^{{\cal P}{\cal P}}_{cd,ef},
\end{array}
\]
hence
\begin{equation}
\label{eq4.1}
R^{{\cal P}{\cal P}}_{ab,cd}(l_{df} * l_{ce}) = (l_{bd} * l_{ac})R^{{\cal
P}{\cal P}}_{cd,ef},\ a,b,e,f = 1,\dots,N,+.
\end{equation}
Setting $l_{ab} = L_{ab}$, $l_{a+} = M_a$, and using $l_{+a}
= 0$, $l_{++} = \e$, $a,b = 1,\dots,N$, and (\ref{eq3.11}), (\ref{eq4.1}) is
equivalent to
\begin{equation}
\label{eq4.2}
R_{ab,cd}(L_{df} * L_{ce}) = (L_{bd} * L_{ac})R_{cd,ef},
\end{equation}
\begin{equation}
\label{eq4.3}
R_{ab,cd}(M_d * L_{ce}) + Z_{ab,c}L_{ce} = (L_{bd} * L_{ac})Z_{cd,e} + L_{be} *
M_a,
\end{equation}
\begin{equation}
\label{eq4.4}
R_{ab,cd}(L_{df} * M_c) - (RZ)_{ab,d}L_{df} = -(L_{bd} * L_{ac})(RZ)_{cd,f} +
M_b * L_{af},
\end{equation}
\begin{equation}
\label{eq4.5}
R_{ab,cd} M_d * M_c + Z_{ab,c}M_c - (RZ)_{ab,d}M_d + s_{ab}\e = (L_{bd} *
L_{ac})s_{cd} + M_b * M_a,
\end{equation}
$a,b,e,f = 1,\dots,N$, where
$s = (R-1)T + m$.  Let us notice that (\ref{eq4.4})
follows from (\ref{eq4.2})--(\ref{eq4.3}).  Moreover, using (\ref{eq1.4}),
(\ref{eq1.3}), $l_{ac}(xy) = l_{ab}(x)l_{bc}(y)$, $l_{ac}(I) = \d_{ac}$,
$a,c = 1,\dots,N,+$, $x,y \in {\cal B}$.  Thus $L_{ac}(I) =
\d_{ac}$, $M_a(I) = 0$,
\begin{equation}
\label{eq4.6}
L_{ac}(xy) = L_{ab}(x)L_{bc}(y),
\end{equation}
\begin{equation}
\label{eq4.7}
M_a(xy) = L_{ab}(x)M_b(y) + M_a(x) \e (y),
\end{equation}
$a,b = 1,\dots,N$, $x,y \in {\cal B}$.  Also $l_{jl}({\cal P}_{ab}) = {\cal
R}({\cal P}_{ab} \otimes {\cal P}_{jl}) = R^{{\cal P}{\cal P}}_{ja,bl}$,
$l_{jl}(w_{AB}) = {\cal R}(w_{AB} \otimes {\cal P}_{jl}) = R^{w{\cal
P}}_{jA,Bl} = (N_w)_{jA,Bl}$ (see $(\ref{eq4.0})$, (\ref{eq1.6'}),
(\ref{eq3.10}), (\ref{eq3.5'})), $a,b,j,l =
1,\dots,N,+$, $A,B = 1,\dots,\dim w$, $w \in \mbox{ Rep } H$.  Therefore
\begin{equation}
\label{eq4.8}
L_{jl}(\Lambda_{ab}) = R_{ja,bl},
\end{equation}
\begin{equation}
\label{eq4.9}
L_{jl}(y_a) = -(RZ)_{ja,l},
\end{equation}
\begin{equation}
\label{eq4.10}
L_{jl}(w_{AB}) = (G_w)_{jA,Bl},
\end{equation}
\begin{equation}
\label{eq4.11}
M_j(\Lambda_{ab}) = Z_{ja,b},
\end{equation}
\begin{equation}
\label{eq4.12}
M_j(y_a) = s_{ja},
\end{equation}
\begin{equation}
\label{eq4.13}
M_j(w_{AB}) = (H_w)_{jA,B},
\end{equation}
$a,b,j,l = 1,\dots,N$, $A,B = 1,\dots,\dim w$, $w \in \mbox{ Rep } H$.  It is
clear that $l_{jl}$ generate a unital subalgebra of ${\cal B}'$ (w.r.t.
convolution $*$).  Endowing it with $\D'$ of Remark~\ref{rem1.4'}, we get a
bialgebra
$U$ with $l$ as its corepresentation.  Adding $l_{ij}^{(m)} = l_{ij} \circ
S^m$, one obtains a Hopf algebra ${\hat U}$ with coinverse $S'(l^{(m)}) =
l^{(m+1)}$ and corepresentations $l^{(2k)}$, $(l^{(2k+1)})^T$,
$k = 0,1,2,\dots$.  Acting $S'{}^m$ on (\ref{eq4.1}), one
obtains
\[
\begin{array}{rll}
R^{{\cal P}{\cal P}}_{ab,cd}(l_{df}^{(2k)} * l_{ce}^{(2k)}) &= &(l_{bd}^{(2k)}
* l_{ac}^{(2k)})R^{{\cal P}{\cal P}}_{cd,ef}, \\
R^{{\cal P}{\cal P}}_{ab,cd}(l^{(2k+1)}_{ce} * l_{df}^{(2k+1)}) &=
&(l_{ac}^{(2k+1)} * l_{bd}^{(2k+1)})R_{cd,ef}^{{\cal P}{\cal P}}.
\end{array}
\]

${\hat U}$ is called enveloping algebra of $({\cal B},\D)$.  It can be
sometimes too small.
It happens e.g. in the classical case (see Remark~\ref{rem3.1})
with ${\cal R}=\e\otimes\e$ when ${\hat
U} = {\cb}\e$. Cf. also \cite{KM}.

Notice that $L_{jl}{|}_{\cal A} = f_{jl}$,
$M_j{|}_{\cal A} = \eta_a$, $l_{{|} {\cal A}} = \rho$.
According to the proof of Theorem~\ref{th1.4}, there exists antihomomorphism
$\th: {\cal B} \to {\cal B}'$ (given by ${\cal R}$) such that
$\th(\Lambda_{jl}) = L_{jl}$, $\th(y_j) = M_j$, $\th(I) = \e$. Therefore
the formulae
$(3.60)$, $(3.46)$ and $(1.12)$ of \cite{INH} yield (\ref{eq4.3}),
(\ref{eq4.5}) (with $s$ replaced by $(R-1)T$)
 and (\ref{eq4.2}) which give
\[
f_{be} * \eta_a = R_{ab,cd}\eta_d * f_{ce} + Z_{ab,c}f_{ce} - (f_{bd} *
f_{ac})Z_{cd,e}
\]
(cf. $(2.18)$ of \cite{INH}), the condition $\tau^{ij} = 0$ and the last
formula before Proposition $3.14$ in \cite{INH}.

Suppose $(\Lambda \otimes \Lambda)k = kw^0$, $n(\Lambda \otimes \Lambda) =
w^0n$ $(w^0 = (I_{\cal B}))$.  Applying $\th$, we get
\begin{equation}
\label{eq4.14}
(L_{bd} * L_{ac})k_{cd} = k_{ab}\e,\ n_{ab}(L_{bd} * L_{ac}) = n_{cd}\e,
\ a,b,c,d
= 1,\dots,N.
\end{equation}

Let us set $X_{lj} = L_{jl} \circ S \in {\hat U}$, $j,l = 1,\dots,N$.

Then
\begin{equation}
\label{eq4.14a}
X_{ik}(xy) = X_{ij}(x)X_{jk}(y),\ X_{ik}(I) = \d_{ik},\ x,y \in {\cal B},
\end{equation}
\begin{equation}
\label{eq4.14b}
X_{ik}(a) = f_{ki}(S(a)),\ a \in {\cal A},\ i,k = 1,\dots,N.
\end{equation}
Using the last equation in the proof of Proposition $3.12$ of \cite{INH},
(\ref{eq4.6}), (\ref{eq3.6}) and (\ref{eq4.9}),
we obtain
\begin{equation}
\label{eq4.14c}
X_{ik}(y_l) = Z_{lk,i}.
\end{equation}
Moreover, (\ref{eq4.14}) yields
\begin{equation}
\label{eq4.14d}
k_{ab}(X_{ac} * X_{bd}) = k_{cd}\e,\ (X_{ac} * X_{bd})n_{cd} = n_{ab}\e.
\end{equation}

As in the proof of Proposition 3.1.2 of \cite{DM}, there exists a unital
homomorphism $X: {\cal B} \to M_{N+1}({\cb})$ such that
\[
X = \left( \begin{array}{cc}
(X_{jl})_{j,l=1}^N & (Y_j)_{j=1}^N \\
0 & \e
\end{array} \right)
\]
for some $Y_j \in {\cal B}'$ satisfying $Y_j(a)=0$, $a\in {\cal A}$,
$Y_j(y_k)=\delta_{jk}$, $j,k=1,\ldots,N$. Setting $X_{j+}=Y_j$, $X_{+j}=0$,
$X_{++}=\e$, $j=1,\ldots,N$,
the commutation relations among $X_{ij}$, $i,j=1,\ldots,N,+$, are the same as
in (3.7) of \cite{DM}, i.e.
\begin{equation}
\label{eq4.18}
(X_{ab} * X_{cd})K_{bd,st} = K_{ac,bd}(X_{bs} * X_{dt}),\ a,c,s,t =
1,\dots,N,+,
\end{equation}
where
\[
K = \left( \begin{array}{cccc}
R^T &0 &0 &0 \\
0 &0 &1\!\!1 &0 \\
0 &1\!\!1 &0 &0 \\
0 &0 &0 &1
\end{array} \right)
\]
(it is also possible to replace $K$ in (\ref{eq4.18}) by $K+n_P$, where
$n\in\mbox{ Mor}(\Lambda\otimes\Lambda,w^0)$ -- see (\ref{eq3.12}),
(\ref{eq4.14d})).

Defining
$X^{(m)}_{ij} = X_{ij} \circ S^m$, $i,j=1,\ldots,N,+$,
$m=0,1,2,\ldots$,  one gets a Hopf algebra
${\hat V}$ generated (as an algebra) by $l_{ij}^{(m)}$ and $X_{ij}^{(m)}$.
Clearly $S'(X^{(m)}) = X^{(m+1)}$; $X^{(2k)}$,
$[X^{(2k+1)}]^T$, $k = 0,1,\dots$, are corepresentations of ${\hat V}$.
Letting
$S'$ act on (\ref{eq4.18}), one obtains
\[
\begin{array}{rll}
(X_{ab}^{(2k)} * X_{cd}^{(2k)})K_{bd,st} &= &K_{ac,bd}(X_{bs}^{(2k)} *
X_{dt}^{(2k)}), \\
(X^{(2k+1)}_{cd} * X_{ab}^{(2k+1)})K_{bd,st} &= &K_{ac,bd}(X_{dt}^{(2k+1)} *
X_{bs}^{(2k+1)}).
\end{array}
\]
${\hat V}$ is called enlarged enveloping algebra of $({\cal B},\D)$.
It would be interesting to find the
commutation relations between $M_i$ and $Y_j$.

\bigskip
\begin{center}
{\bf Acknowledgment}
\end{center}

\medskip
The author is grateful to Professor W. Arveson and other faculty members for
their kind hospitality at UC Berkeley.


\begin{thebibliography}{99}

\bibitem[1]{Dr} Drinfeld, V. G., Quantum groups, in {\em Proceedings
ICM-1986, Berkeley}, (1987), 798--820.

\bibitem[2]{Ha} Hayashi, T., Quantum groups and quantum determinants,
{\em Journal of Algebra}, {\bf 152:1} (1992), 146--165.

\bibitem[3]{KM} Kosi\'nski, P. and Ma\'slanka, P., The duality between
$\kappa$-Poincar\'e algebra and
$\kappa$-Poincar\'e group, {\em hep-th 9411033}.

\bibitem[4]{LT} Larson, R. G. and Towber, J., Two dual classes of bialgebras
related to the concepts of ``quantum group'' and ``quantum Lie algebra'',
{\em Comm. Algebra}, {\bf 19} (1991), 3295--3345.

\bibitem[5]{Lyu} Lyubashenko, V. V., Hopf algebras and vector symmetries,
{\em Russian Math. Surveys}, {\bf 41:5} (1986), 153--4. Russian original: {\em
Uspekhi Mat. Nauk}, {\bf 41:5} (1986), 185--186.

\bibitem[6]{Ma} Majid, S., Quasitriangular Hopf algebras and Yang--Baxter
equations, {\em Int. Journal of Modern Phys. A}, {\bf 5:1} (1990), 1--91.

\bibitem[7]{CQG} Podle\'s, P., Complex quantum groups and their real
representations, {\em Publ. RIMS, Kyoto Univ.}, {\bf 28} (1992), 709--745.

\bibitem[8]{DM} Podle\'s, P., Solutions of Klein--Gordon and Dirac
equations on quantum
Minkowski spaces, q-alg 9510019, UC Berkeley preprint, PAM-654.

\bibitem[9]{L} Podle\'s, P. and Woronowicz, S. L., Quantum deformation of
Lorentz group, {\em Commun. Math. Phys.}, {\bf 130} (1990), 381--431.

\bibitem[10]{INH} Podle\'s, P. and Woronowicz, S. L., On the structure of
inhomogeneous quantum groups, hep-th 9412058, UC Berkeley preprint,
PAM-631.

\bibitem[11]{POI} Podle\'s, P. and Woronowicz, S. L., On the classification of
quantum Poincar\'e groups, hep-th 9412059, UC Berkeley preprint, PAM-632.

\bibitem[12]{RTF} Reshetikhin, N. Yu., Takhtadzyan, L. A. and Faddeev, L. D.,
Quantization of Lie groups and Lie algebras, {\em Leningrad Math. J.}, {\bf
1:1} (1990), 193--225.  Russian original:  {\em Algebra i analiz}, {\bf 1:1}
(1989), 178--206.

\bibitem[13]{UEA} Schlieker, M., Weich, W. and Weixler, R., Inhomogeneous
quantum groups and their
quantized universal enveloping algebras, {\em Lett. Math.
Phys.}, {\bf 27} (1993), 217--222.

\bibitem[14]{Ta} Takeuchi, M., Finite dimensional representations of
the quantum Lorentz group, {\em Commun. Math. Phys.}, {\bf 144} (1992),
557--580.

\bibitem[15]{W1} Woronowicz, S. L., Twisted $SU(2)$ group.  An example of a
non-commutative differential calculus, {\em Publ. RIMS, Kyoto Univ.}, {\bf 23}
(1987), 117--181.

\bibitem[16]{TK} Woronowicz, S. L., Tannaka--Krein duality for compact matrix
pseudogroups.  Twisted $SU(N)$ groups, {\em Inv. Math.}, {\bf 93} (1988),
35--76.

\bibitem[17]{WP} Woronowicz, S.L., private information.

\bibitem[18]{WZ} Woronowicz, S. L. and Zakrzewski, S., Quantum deformations of
the Lorentz group. The Hopf ${}^*$-algebra level, {\em Comp. Math.}, {\bf 90}
(1994), 211--243.

\end{thebibliography}
\end{document}